\documentclass[12pt]{amsart}

\usepackage[centertags]{amsmath}

\usepackage{amsfonts}
\usepackage{amssymb}
\def\RR{\mathbb R}
\def\QQ{\mathbb Q}
\def\FF{\mathbb F}

\def\into{\longrightarrow}

\def\res{\hspace{-4pt} \upharpoonright \hspace{-3pt}}
\def\harp{\res}

\DeclareMathOperator{\supp}{supp}

\DeclareMathOperator{\dom}{dom}

\title{Chromatic Numbers of \\ algebraic hypergraphs}

\author{James H. Schmerl}

\date{\today}

\begin{document}

\begin{abstract}
Given a polynomial $p(x_0,x_1, \ldots, x_{k-1})$ over the reals $\RR$, where each 
$x_i$ is an $n$-tuple of variables, we form its zero $k$-hypergraph $H = (\RR^n,E)$,
where the set $E$ of edges consists of all $k$-element sets $\{a_0,a_1, \ldots, a_{k-1}\} \subseteq \RR^n$ such that $p(a_0,a_1, \ldots, a_{k-1}) = 0$.  Such hypergraphs are precisely the {\it algebraic} hypergraphs. We say (as in \cite{avoid}) that $p(x_0,x_1,$ $ \ldots, x_{k-1})$ 
is {\it avoidable}  if the  chromatic number $\chi(H)$ of its zero hypergraph $H$ is countable, and it is $\kappa$-{\it avoidable} if $\chi(H) \leq \kappa$. Avoidable polynomials were 
completely characterized in \cite{avoid}. 
For any infinite $\kappa$, we  characterize the $\kappa$-avoidable 
algebraic hypergraphs. 
Other results about algebraic hypergraphs and their chromatic numbers are also proved.
\end{abstract}

\maketitle

 A polynomial $p(x_0,x_1, \ldots, x_{k-1})$ over the reals $\RR$ is $(k,n)$-{\it ary} if each $x_i$ 
is an $n$-tuple of variables.  Following \cite{avoid}, we say  that a $(k,n)$-ary polynomial $p(x_0,x_1, \ldots, x_{k-1})$ is  {\em avoidable} if the points of $\RR^n$ can be colored with countably many colors 
 such that whenever $a_0,a_1, \ldots, a_{k-1} \in \RR^n$ 
are distinct and $p(a_0,a_1, \ldots, a_{k-1}) = 0$,  
then there are  $i < j < k$ such that the points $a_i,a_j$ are   differently  colored.  The avoidable polynomials were characterized in \cite{avoid}. 
The prototypical examples  of avoidable polynomials are the $(3,n)$-ary polynomials $\|x-y\|^2 - \|y-z\|^2$, where $n \geq 2$,
which were shown in \cite{96} to be avoidable after various partial results had been obtained in \cite{t1},\cite{t3},\cite{t2},\cite{t4},\cite{t5}. A consequence with a more geometric flavor 
might be stated succinctly as {\em 
 the set of isosceles triangles 
is avoidable.} Another example is the 
 $(3,2)$-ary polynomial $\|x-y\|^2 + \|y-z\|^2 - \|x-z\|^2$, which  is avoidable iff the Continuum Hypothesis  ({\sf CH}) is true.  Fox \cite[Coro.\@ 1]{fox} showed that the homogeneous linear $(k+3,1)$-ary polynomial $x_0 + x_1 + \cdots + x_k - x_{k+1} - kx_{k+2}$ is avoidable iff 
$2^{{\aleph}_0} \leq \aleph_k$. Other examples can be found in  \cite{avoid}.

Observe that the  polynomial $p(x_0,x_1, \ldots, x_{k-1})$ is avoidable iff a certain 
 hypergraph (which we will refer to as its zero hypergraph) has a countable chromatic number.  There was no attempt in  \cite{avoid} to say anything   additional about the chromatic numbers of the zero hypergraphs of unavoidable polynomials. Building on \cite{avoid}, we will rectify this omission in \S2 by  saying exactly what these chromatic numbers are. 

Recall that $H = (V,E)$ is  a $k$-{\it hypergraph} (or a $k$-{\it uniform} hypergraph, as it is more usually referred to) if $k < \omega$, $V$  is a nonempty set and 
$E$  is a  set of $k$-element subsets of $V$.  For such an $H$, $V$ is its set of {\it vertices} and $E$ its set of {\it edges} (or {\it hyperedges}). All hypergraphs considered here are assumed to be $k$-hypergraphs for some $k < \omega$. We will usually assume that  $k \geq 2$. 
A {\it graph} is just a $2$-hypergraph. A  function $\varphi : V \into C$ 
is a $\kappa$-{\it coloring} of $H$ if $|C| \leq \kappa$, and it is {\it proper coloring} of $H$ provided that $\varphi$ is not constant on any edge of $H$. If there is a proper $\kappa$-coloring of $H$, then $H$ is $\kappa$-{\it colorable}. 
The {\it chromatic number} $\chi(H)$ of $H$ is the least  (possibly infinite) cardinal $\kappa$ 
such that $H$  is $\kappa$-colorable. 

A subset $X \subseteq \RR^n$ is {\it algebraic} if it is the zero-set of a polynomial 
(or, equivalently, of a finite set of polynomials) over $\RR$. 
If $p(x_0,x_1, \ldots, x_{k-1})$ is a $(k,n)$-ary polynomial and $H = (V,E)$ is a $k$-hypergraph, 
then we say that $H$ is the {\it zero} hypergraph of $p(x_0,x_1, \ldots, x_{k-1})$ if 
$V = \RR^n$ and $E = \{\{a_0,a_1, \ldots, a_{k-1}\} \subseteq \RR^n : |\{a_0,a_1, \ldots, a_{k-1}\}| = k$ 
and $p(a_0,a_1, \ldots,a_{k-1}) = 0\}$. A hypergraph is {\it algebraic} if it is the zero hypergraph of  some polynomial over $\RR$. A polynomial $p(x_0,x_1, \ldots, x_{k-1})$ is $\kappa$-{\it avoidable}
if its zero hypergraph  has a proper $\kappa$-coloring. 

A famous  example of an algebraic graph  is the {\it unit-distance graph} ${\bf X}_2(\{1\}) = (\RR^2,E)$, whose vertices are the points in the plane and whose edges are those pairs of points at a distance 1 from each other. The notoriously  obstinate Hadwiger-Nelson problem is to determine  $\chi({\bf X}_2(\{1\}))$.   For an entertaining account 
of this problem and its history, see \cite[Chaps.\@ 2 \& 3]{soi}. All that is known 
about the exact value of $\chi({\bf X}_2(\{1\})$ are the relatively easy bounds of 
$4 \leq \chi({\bf X}_2(\{1\})) \leq 7$ that were established soon after the problem was proposed 
in the early 1950's. For every $n < \omega$, the unit-distance graph on $\RR^n$,
denoted by ${\bf X}_n(\{1\})$, has a finite chromatic number.   
The  determination of   finite chromatic numbers  is  essentially a finite problem 
since, by the De~Bruijn-Erd\H{o}s Theorem, if $H$ is a hypergraph and $n < \omega$, then 
$H$ is $n$-colorable  iff every finite subhypergraph of $H$ is $n$-colorable. However, the focus of this paper  will  be on infinite chromatic numbers.

There are three numbered sections following this introduction. The first, \S0,  contains some preliminary material that can be skipped by most readers. The  simplest of all 
algebraic hypergraphs -- those obtained from templates -- are defined and studied in \S1
where their chromatic numbers are determined exactly. The principal results, linking chromatic numbers of algebraic hypergraphs to the chromatic numbers of the 
hypergraphs defined in \S1,  appears in \S2.  

\bigskip


{\bf \S0.\@ Preliminaries.}  This section contains some  definitions that should be familiar.

If $\alpha$ is an (infinite or finite) ordinal, then $\alpha$ is the set of its predecessors. 
As usual, $\omega$ is the least infinite ordinal. Thus, $\omega$ is the set of natural numbers.
For any set $X$ and $n < \omega$, $X^n$ is the set of $n$-tuples from $X$.
If $x \in X^n$, then we will often understand that $x = \langle x_0,x_1, \ldots, x_{n-1} \rangle$. 
If $x$ is an $n$-tuple and $a$ is some element, then $xa$ is that $(n+1)$-tuple extending $x$ 
such that $(xa)_n = a$.

If $n,k < \omega$, then, in the appropriate context, $n^k$ will   be $\{0,1, \ldots,n-1\}^k$, which is the set of $k$-tuples of elements of $n$. If $\alpha$ is an infinite ordinal or cardinal, then $\alpha^k$ will always be the set of $k$-tuples of ordinals that are less than $\alpha$.
        
        If $X$ is any set and $n < \omega$, then $[X]^n$ is the set of all $n$-element subsets of $X$. 
        Thus, $(V,E)$ is a $k$-hypergraph iff $E \subseteq [V]^k$. If $X$ is linearly ordered by $<$ 
        (for example, if $X$ is an ordinal) and we write that $\{x_0,x_1, \ldots,x_{n-1}\}_< \subseteq X$, then it is to be understood that 
     $x_0 < x_1 < \cdots < x_{n-1}$. We will use the notational expedient of letting      
      $x = \langle x_0,x_1, \ldots, x_{n-1}\rangle$ whenever we have already agreed that $\{x_0,x_1, \ldots, x_{n-1}\}_< \in [X]^n$.

If $\kappa$ is a (typically,  infinite) cardinal, then $\kappa^+$ is the successor cardinal of $\kappa$. If $n < \omega$, 
then $\kappa^{+n}$ is defined recursively by $\kappa^{+0} = \kappa$ and $\kappa^{+(n+1)} = (\kappa^{+n})^+$. Thus, $(\aleph_\alpha)^{+n} = \aleph_{\alpha+n}$. We could define $\kappa^{+\alpha}$ for any ordinal $\alpha$, but the only infinite ordinal we need is $\alpha = \omega$, in which case $\kappa^{+\omega} = \bigcup\{\kappa^{+n} : n < \omega\}$.
     
     Let  $\RR$ be the set of reals and $\widetilde \RR = (\RR,+,\times,0,1,\leq)$  be the ordered field of the  real numbers. 
Let ${\mathcal L}_{OF} = \{+,\times,0,1,\leq\}$ be the first-order language appropriate for ordered fields. If $D \subseteq \RR$, let ${\mathcal L}_{OF}(D)$ be ${\mathcal L}_{OF}$ augmented with names for the elements 
of $D$. If $n < \omega$ and $A \subseteq \RR^n$, then $A$ is $D$-{\it definable} if it is definable in $\widetilde \RR$ by a first-order ${\mathcal L}_{OF}(D)$-formula, and it is   {\it semialgebraic} iff  it is $\RR$-definable.  The algebraic sets  were defined in the introduction.

  Suppose that $H_1 = (V_1,E_1)$ and $H_2 = (V_2,E_2)$ are $k$-hypergraphs. 
As usual, $H_1$ is a {\em subhypergraph} of $H_2$ if $V_1 \subseteq V_2$ and $E_1 \subseteq E_2$. If $f : V_1 \into V_2$, then $f$ {\it embeds} $H_1$ into $H_2$ if $f$ is an isomorphism 
onto a subhypergraph of $H_2$. We say that $H_2$ {\em contains} an $H_1$ if some $f$ embeds 
$H_1$ into $H_2$.  Obviously, if $H_2$ contains an $H_1$, then $\chi(H_2) \geq \chi(H_1)$.

\bigskip


{\bf \S1.\@ Templates and Their Hypergraphs.} If $1 \leq d < \omega$ and $2  \leq  k < \omega$, then $P$ is a $d$-{\it dimensional} $k$-{\it template} if $P$ is a set of $d$-tuples and $|P| = k$. Two $d$-dimensional $k$-templates $P,Q$ are {\it isomorphic} if there is a bijection $f : P \into Q$ such that whenever $x,y \in P$  and $i < d$, 
then $x_i = y_i$ iff  $f(x)_i = f(y)_i$. 
 Obviously, for fixed  $d$ and $k$,  there are only finitely many, non-isomorphic $d$-dimensional $k$-templates. If both $P$ and $Q$ are $d$-dimensional templates, then we say that $Q$ is a {\it homomorphic image} of~$P$ if there is a surjective  function $f : P \into Q$ such that   whenever  $x,y \in P$, $i < d$ 
and $x_i = y_i$, then $f(x)_i = f(y)_i$.

If  $X = X_0 \times X_1 \times \cdots \times X_{d-1}$ and $P$ is  a $d$-dimensional $k$-template, then its {\it template hypergraph} $L(X,P)$ on $X$ is the $k$-hypergraph whose set of vertices is $X$ and whose edges are those $k$-templates $Q \subseteq X$ that are homomorphic images of~$P$. 
If $P$ is a $d$-dimensional $k$-template, then $L(\RR^d,P)$ is an algebraic $k$-hypergraph.
For example, if $P$ is the $2$-dimensional $3$-template $\{\langle 0,0 \rangle, \langle 0,1 \rangle, \langle 1,1 \rangle\}$, then $L(\RR^2,P)$ is the 3-hypergraph that is the zero hypergraph of  the $(3,2)$-ary polynomial 
$p(x,y,z) = (x_0-y_0)^2 + (y_1 - z_1)^2$. 

In this section, the chromatic numbers of the various $L(\RR^d,P)$ will be determined.

Let $P$ be a $d$-dimensional $k$-template. We say that a  subset $I \subseteq d$ is a {\em distinguisher} for 
$P$ if whenever $x,y \in P$ are distinct, then $x_i \neq y_i$ for some $i \in I$. We then define 
$e(P)$ to be the least $e$ that is the cardinality of a distinguisher. Obviously, $e(P) \leq d$ since $d$ 
itself is a distinguisher. In addition, $1 \leq e(P) \leq k-1$. The lower bound is trivial since we are assuming that $k \geq 2$.  The upper bound 
$e(P) \leq k-1$ is proved by induction on $k$. It is obvious for $k = 2$ (or even $k = 1$). 
Now suppose that $k \geq 3$ and we have proved the inequality for all smaller~$k$. 
Let $x \in P$ and $Q = P \backslash \{x\}$. By the inductive hypothesis, there is a distinguisher $I$ for $Q$ such that $|I| \leq k-2$. If $I$ is a distinguisher for $P$ we are done. Otherwise, there is a unique $y \in Q$ such that $x \harp I = y \harp I$. Let $i < d$ be such that $x_i \neq y_i$, and then $I \cup \{i\}$ is a distinguisher for $P$ and has cardinality at most $k-1$. 
If $Q$ is a $d$-dimensional $k$-template that is a homomorphic image of $P$, then $e(Q) \geq e(P)$ 
since any distinguisher for $Q$ is also one for $P$.

The next theorem characterizing 
each $\chi(L(\RR^d,P))$  is the main result of this section.

\bigskip

{\sc Theorem 1.1}: {\em Suppose that  $P$ is a $d$-dimensional $k$-template. Then $\chi(L(\RR^d,P))$ is the least $\kappa$ such that 
$\kappa^{+(e(P)-1)} \geq 2^{\aleph_0}$.}

\bigskip

The proof of Theorem~1.1 will be given after several supporting lemmas and corollaries.
\bigskip

{\sc Lemma 1.2:} {\em Suppose that  $\kappa \geq \aleph_0$,  
$P$ is a  $d$-dimensional $k$-template, and $X_0,X_1, \ldots, X_{d-1}$ 
are sets such that $|X_i| \geq \kappa^{+i}$ for each $i < d$. 
Let $X = X_0 \times X_1 \times \cdots \times X_{d-1}$.
 Then 
$\chi(L(X,P)) \geq \kappa$.}

\bigskip

{\it Proof}.  The proof is by induction on the cardinal $\kappa$. By the multipartite Ramsey theorem 
(which is midway in strength between Theorems~1 and 5 of \cite[Chap.\@ 5.1]{grs}), the lemma is true when $\kappa = \aleph_0$. 
(In fact, all that is required of the $X_i$'s is that they be infinite.) The instance of the lemma when $\kappa$ is an uncountable limit cardinal will follow from all instances for smaller infinite cardinals, so we can assume that $\kappa = \lambda^+$. Then the lemma is  essentially the theorem of Erd\H{o}s \& Hajnal (see \cite[Lemma~1.1]{avoid})
which asserts: If $|X_i| =  \lambda^{+(i+1)}$ for all $i < d$ and $f : X \into \lambda$, 
then for every $t < \omega$, there are $B_i \subseteq X_i$ for each $i < d$ such that 
$|B_0| = |B_1| = \cdots = |B_{d-1}| = t$ and $f$ is constant on $B_0 \times B_1 \times \cdots \times B_{d-1}$. Consider arbitrary $f : X \into \lambda$ and choose $t = dk$. We have chosen $t$ to be large enough so that  there are 
$b_0,b_1, \ldots, b_{k-1} \in B_0 \times B_1 \times \cdots \times B_{d-1}$ such that 
$Q = \{ b_0,b_1, \ldots, b_{k-1}\}$  is a $k$-template that is isomorphic to $P$. Hence,  $Q$ is an edge of 
$L(X,P)$ and 
$f(b_0) = f(b_1) = \cdots = f(b_{k-1})$, so $f$ is not a proper coloring of $L(X,P)$, 
thereby showing that {\mbox{$\chi(L(X,P)) \geq \lambda^+$}.
\qed

\bigskip

{\sc Corollary 1.3:} {\em If $P$ is a  $d$-dimensional template,
then 

\vspace{11pt}
 
\hspace{90pt}
$\chi(L(\RR^d,P))^{+(d-1)} \geq 2^{\aleph_0}$.}

\bigskip

{\it Proof}. Let $\lambda = \chi(L(\RR^d,P))$ and, for a contradiction, suppose that 
$\lambda^{+(d-1)} < 2^{\aleph_0}$. Then, $(\lambda^+)^{+(d-1)} \leq 2^{\aleph_0}$, so Lemma~1.2 implies (with $\kappa = \lambda^+$) that $\chi(L(\RR^d,P)) \geq \lambda^+$, a contradiction. \qed

 \bigskip

We will say that a $d$-dimensional template $P$ is {\it simple} 
if for every $m < d$, there are $x,y \in P$ such that for all $i <d$, $ x_i = y_i$ iff $i \neq m$. 
We easily see, by induction on $d$, that if $P$ is a simple $d$-dimensional $k$-template, then $k \geq d+1$. 
Obviously, for fixed  $k < \omega$,  there are only finitely many non-isomorphic simple $k$-templates. An example of a simple  $d$-dimensional $(d+1)$-template is 
$\{ x_0,x_1, \ldots, x_d \}$, where each $x_i$ is a $d$-tuple of all $0$'s except that $x_{ii}=1$ if $i < d$. 
To get a simple $d$-dimensional $k$-template with $k > d+1$ just add $k-d-1$ new 
$d$-tuples to  the  aforementioned example.

\bigskip

{\sc Lemma 1.4:} {\em Suppose that $P$ is a  simple $d$-dimensional $k$-template and $\kappa \geq \aleph_0$.
If $|X| \leq \kappa^{+(d-1)}$, 
then   $\chi(L(X^d,P)) \leq  \kappa$.}

\bigskip

{\it Proof}. We make use of a result of Komj\'{a}th (see \cite[Lemma~1.3]{avoid}) that asserts:
If $|X| \leq  \kappa^{+(d-1)}$, then there are functions 
$G : X^d \into \kappa$ and $j : X^d \into d$ such that whenever $a,b \in X^d$, $G(a) = G(b)$,  
$j(a) = j(b)$  and 
$a_i = b_i$ for $j(a) \neq i \leq d$, then $a = b$. Having such $G$ and $j$, we let 
$f$ be the function  on $X^d$ such that $f(a) = \langle G(a),j(a) \rangle$. 
Clearly, $f$ is a $\kappa$-coloring of $L(X^d,P)$. To see that $f$   is  proper, 
suppose $\{x_0,x_1, \ldots, x_{k-1}\}$ 
is an edge and 
$f(x_0) = f(x_1) = \cdots = f(x_{k-1}) = \langle \alpha, m \rangle$. Let $r < s < k$ be such 
for every $i < d$, $x_{r,i} = x_{s,i}$ iff $i \neq m$. Then, $x_r = x_s$, a contradiction. 
\qed

\bigskip

{\sc Corollary 1.5:} {\em Let $P \subseteq \RR^d$ be a   $k$-template.}

(a) {\em If $P$ is simple, 
then $\chi(L(\RR^d,P))$  is the least $\kappa$ such that $\kappa^{+(d-1)} \geq 2^{\aleph_0}$.}

(b) {\em Suppose that $I \subseteq d$ is a distinguisher and that $|I| = e(P)$. 
Let $Q = \{x \in \RR^d : x \harp I = y \harp I$ for some $y \in P$ and $x_i = 0$ for $i \in d \backslash I\}$.
Then $\chi(L(\RR^d,Q))$  is the least $\kappa$ such that $\kappa^{+(e(P)-1)} \geq 2^{\aleph_0}$.}

\bigskip

{\it Proof.} (a) is immediate from Corollary~1.3 and Lemma~1.4. 

To prove (b), let $e = e(P)$ and let $\kappa$ be the least such that $\kappa^{+(e-1)} \geq 2^{\aleph_0}$. Without loss of generality, suppose that $I = e$. Let $R = \{x \harp e : x \in P\}$. Then $R$ is 
a simple $e$-dimensional $k$-template. By (a), $\chi(L(\RR^e,R)) = \kappa$. 
We will show that $\chi(L(\RR^d,Q)) = \chi(L(\RR^e,R))$. 

 Consider any proper coloring $\varphi$ of $L(\RR^e,R))$. Let $\theta$ be the unique coloring 
of $L(\RR^d,Q)$ such that $\theta(x) = \varphi(x \harp e)$.  
Then, $\theta$ is a proper coloring of $L(\RR^d,Q)$ so that $\chi(L(\RR^d,Q)) \leq \chi(L(\RR^e,R))$.
For the reverse  inequality, consider any proper coloring $\theta$ of $L(\RR^d,Q)$. 
Then let $\varphi$ be the unique coloring of $L(\RR^e,R))$ such that 
$\varphi(x) = \theta(y)$, where $y \harp e = x$ and $y_i = 0$ for $e \leq i < d$. 
Then, $\varphi$ is a proper coloring of $L(\RR^e,R)$ so that $\chi(L(\RR^d,Q)) \geq \chi(L(\RR^e,R))$.
\qed

\bigskip

The next corollary is for later use (in proving Corollary~2.3).

\bigskip

{\sc Corollary 1.6}: {\em Suppose that $2 \leq k < \omega$ and $\kappa \leq 2^{\aleph_0} \leq \kappa^{+(k-2)}$. Then there is a $d$-dimensional $k$-template $P$ such that $\chi(L(\RR^d,P)) = \kappa$.}

\bigskip

{\it Proof}. 
 Let  $d$ be such that 
$1 \leq d <k$ and 
$\kappa^{+(d-1)} = 2^{\aleph_0}$, and let $P$ be  a simple $d$-dimensional 
$k$-template}. Apply Corollary~1.5(a). 
\qed

\bigskip

{\it Proof of Theorem~1.1}. Let $P,d,k$ be as in the Theorem. Let $e = e(P)$ and let $\kappa$ be the least such that $\kappa^{+(e-1)} \geq 2^{\aleph_0}$. 

We first show that $\kappa \leq \chi(L(\RR^d,P))$. Let $I \subseteq d$ be a distinguisher  for $P$  such that  $|I| = e$.  Without loss of generality, assume that $I = e \leq d$. Let $Q$ be the homomorphic image of $P$ consisting of those $x \in \RR^d$ such that for some $y \in P$, $x \harp e = y \harp e$ and $x_i = 0$ for $e \leq i < d$. By Corollary~1.5(b), $\chi(L(\RR^d, Q)) = \kappa$. But $L(\RR^d, Q)$ is embeddable into $L(\RR^d,P)$, so $\chi(L(\RR^d,P)) \geq \kappa$.

Next, we prove that $\kappa \geq \chi(L(\RR^d,P))$ by exhibiting a proper $\kappa$-coloring of 
$L(\RR^d,P)$.
Let $Q_0,Q_1, \ldots, Q_m \subseteq \RR^d$ be all (up to isomorphism)  of the  homomorphic images of $P$ 
that are $d$-dimensional $k$-templates.  For each $j \leq m$, let $I_j \subseteq d$ be a distinguisher for $Q_j$ such that $|I_j| = e_j = e(Q_j)$. 
Each $e_j \geq e$, so if $\kappa_j$ is the least such that $\kappa_j^{+(e_j -1)} \geq 2^{\aleph_0}$, 
then $\kappa_j \leq \kappa$. 
By Corollary~1.5(b),  $\chi(L(\RR^d,Q_j)) = \kappa_j \leq \kappa$. Let $\varphi_j$ be a proper $\kappa$-coloring of $L(\RR^d,Q_j)$.

We are now prepared to obtain a proper $\kappa$-coloring $\varphi$ of $L(\RR^d,P)$. For each $x \in \RR^d$,
let
$$        
\varphi(x) = \langle \varphi_j(x)) : j \leq m \rangle.
$$
Clearly, $\varphi$ is a $\kappa$-coloring of $L(\RR^d.P)$. To see that it is proper, consider an edge $Q $ of 
$L(\RR^d,P)$. Then $Q$ is a homomorphic image of $P$, so we can let $j \leq m$ be such that 
$L(\RR^d,Q) = L(\RR^d,Q_j)$. There are distinct $x,y \in Q$ such that $\varphi_j(x) \neq \varphi_j(y)$. But then 
 $\varphi(x) \neq \varphi(y)$. \qed

\bigskip

We end this section with two results about embedding some template hypergraphs into others.

\bigskip

 {\sc Lemma 1.9}:  {\em Let $P$ be a $d$-dimensional $k$-template and let $e = e(P)$. There is  a simple $e$-dimensional $k$-template $Q$ such that for any set $X$, $L(X^e,Q)$ is embeddable into 
$L(X^d,P)$.}

\bigskip

{\it Proof}. Let $I \subseteq d$ be a distinguisher for $P$ such that $|I| = e$. Without loss of generality,
assume that $I = e$. Let $Q = \{x \harp e : x \in P\}$. Clearly, $Q$ is a simple $e$-dimensional 
$k$-template. 

To see that $Q$ is as required, consider any nonempty set $X$ and let $a \in X$. 
Let $f : X^e \into X^d$ 
be such that if $x \in X^e$, then $f(x) = y \in X^d$, where $x = y \harp e$ and 
$y_i = a$ for $e \leq i < d$. 
It is easily checked that $f$ embeds $L(X^e,Q)$  into 
$L(X^d,P)$. \qed

\bigskip

{\sc Lemma 1.10}: {\em Let $P$ be a $d$-dimensional $k$-template and let $m \geq d$. Then there is 
an $m$-dimensional $k$-template $Q$ such that  $e(Q) = e(P)$, $P = \{x \harp d : x \in Q\}$ and for any infinite $X$,
$L(X^m,Q)$ is embeddable into $L(X^d,P)$.}

\bigskip

{\it Proof}. When $m = d$, let $Q = P$. For $m > d$, by an inductive proof it suffices to 
let $m=d+1$. We fix $d \geq 1$ and then prove this case by induction on $k$. 

 Let $X$ be an infinite set. Partition $X$ into $|X|$ subsets each of cardinality $|X|$. Index these sets by elements of $X$.
Thus, we have a partition $\{X_a : a \in X\}$ of $X$ where each $|X_a| = |X|$. 
For each $a \in X$, let $f_a : X \into X_a$ be a bijection. We now define $g_{X} : X^{d+1} \into X^d$ so such that if $a_0,a_1, \ldots, a_d \in X$, then 
$$g_{X}(a_0,a_1, \ldots,a_d) =  \langle f_{a_d}(a_0), f_{a_d}(a_1), \ldots, f_{a_d}(a_{d-1}) \rangle.$$ 
Clearly, $g_{X}$ is an injection.

We will need another  definition. We say that the $d$-dimensional $k$-template $P$ is {\it connected} if, whenever $P$ is  partitioned into two nonempty subsets $P_0,P_1$, then there are  
 $x \in P_0$, $y \in P_1$ and $j < d$ such that $x_j = y_j$. If $P$ is connected, then so is every 
 $k$-template that is a homomorphic image of $P$.

For $2 \leq k < \omega$, consider the following statement.

\begin{quote}

\hspace{-33pt} $(\bigstar_k)$ For every $d$-dimensional $k$-template $P$ there is a $(d+1)$-dimensional $k$-template $Q$ 
such that  $P = \{x \harp d : x \in Q\}$ and for any infinite $X$, $g_{X}$ embeds 
$L(X^{d+1},Q) \into L(X^d,P)$. Also, $e(Q) = e(P)$; moreover, if $I \subseteq d+1$ is a distinguisher of $Q$, then, for each $j < d$,  
$(I \cap d) \cup \{j\}$ is a distinguisher for $P$.. 

\end{quote}

We will prove  $({\bigstar_k})$ for all $k$ by induction on $k$.

\smallskip

The basis step:  $k = 2$. Let $P = \{x,y\}$, where $x,y$ are distinct $d$-tuples 
of $0$'s and $1$'s. Let $X$ be an infinite set. We will write $g$ for $g_{X}$. There are two cases to consider, depending on whether or not $P$ is connected.

\smallskip

{\em $P$ is connected.}  Thus, there is $j < d$ such that $x_j = y_j$. 
Let $Q = \{x0,y0\}$. Clearly, $Q$ is $(d+1)$-dimensional $2$-template such that $e(Q) = e(P) = 1$ and  $P = \{x \harp d : x \in Q\}$. The ``moreover'' part of $(\bigstar_2)$ is obvious. 
The argument that $g$ embeds $L(X^{d+1},Q)$ into $L(X^d,P)$ is straightforward. 
Consider $T = \{x,y\} \in [X^{d+1}]^2$, which is either an edge of $L(X^{d+1},Q)$ or is not.

Suppose that $T$ is an edge of $L(X^{d+1},Q)$. Then $x_d = y_d = a$ for some $a \in X$, and 
$\{x \harp d, y \harp d\}$ is an edge of $L(X^d,P)$. Then 
$g(x) = \langle f_a(x_0),f_a(x_1), \ldots, f_a(x_{d-1}) \rangle$ and 
$g(y) =  \langle f_a(y_0),f_a(y_1), \ldots, f_a(y_{d-1}) \rangle$. Since $f_a : X \into X_a$ is a bijection, 
$\{g(x), g(y)\}$ is an edge of $L(X_a^d,P)$ and, hence, also of $L(X^d,P)$. 

Suppose that $T$ is not an edge of $L(X^{d+1},Q)$. This is due to either $x_d \neq y_d$ or else 
$x_d = y_d$ and $\{x \harp d, y \harp d\}$ is not an edge of $L(X^d,P)$. In the first case, we get that 
$\{g(x),g(y)\}$ is not connected, so it is not an edge of $L(X^d,P)$. In the second case, we have $x_d = y_d = a$ for some $a \in X$. 
Again, $f_a : X \into X_a$ is a bijection, implying that $\{g(x),g(y)\}$ is not an edge of $L(X_a^d,P)$ and, hence, also not an edge of $L(X^d,P)$. 

\smallskip

{\em $P$ is not connected.}
We then let $Q = \{x0,y1\}$.   Clearly, $Q$ is $(d+1)$-dimensional $2$-template such that $e(Q) = e(P) = 1$ and  $P = \{x \harp d : x \in Q\}$. The ``moreover'' part of $(\bigstar_2)$ is clear. Let $X$ be infinite. Observe, in this case, that $L(X^d,P)$ and $L(X^{d+1},Q)$ are both 
complete graphs, so clearly $g$ is an embedding.

\smallskip

The inductive step.  Consider $k \geq 3$ and assume that $(\bigstar_\ell)$ is true whenever $2 \leq \ell < k$. We will prove $(\bigstar_k)$. The proof is similar to the proof of the basis step. 
Let $P$ be a $d$-dimensional $k$-template.  Let $X$ be an infinite set. 
We will write $g$ for $g_{X}$.

We consider two cases, depending on whether or not $P$ is connected. 

\smallskip

 {\em $P$ is connected.} Let $Q = \{x0 : x \in P\}$. Clearly, $Q$ is a $(d+1)$-dimensional $k$-template 
 such that $e(Q) = e(P)$ and $P = \{x \harp d : x \in Q\}$. The ``moreover'' part of $(\bigstar_2)$ is obvious.
 We will show that $g_X$ is an embedding of  $L(X^{d+1},Q)$  into $L(X^d,P)$.
 Consider $T = \{t_0,t_1, \ldots, t_{k-1}\} \in [X^{d+1}]^k$. 
 Either $T$ is an edge of $L(X^{d+1},Q)$ or not.

Suppose that $T$ is an edge of $L(X^{d+1},Q)$. 
There is $a \in X$ such that $t_{0,d} = t_{1,d} = \cdots = t_{k-1,d} = a$, and $\{x \harp d : x \in Q\}$ is an edge of $L(X^d,P)$. 
For each $i < k$, $g(t_i) = \langle f_a(t_{i,0}), f_a(t_{i,1}), \ldots, f_a(t_{i,d-1}) \rangle$.
Since $f_a : X \into X_a$ is a bijection, then  $g[T]$ is an edge of $L(X_a^d,P)$,
so it also is an edge of $L(X^d,P)$.

Suppose that $T$ is not an edge of $L(X^{d+1},Q)$. This failure  is due to one of two reasons: either there are $x,y \in T$ such that $x_d \neq y_d$ or else $x_d = y_d$ for all $x,y \in T$ and $\{x \harp d : x \in T\}$ is not an edge 
of $L(X^d,P)$. 

In the first case, let $x,y \in T$ be such that $x_d \neq y_d$. Let $T_0 = \{ z \in T :  z_d = x_d\}$ and $T_1 = T \backslash T_0$. Thus, $T_0,T_1 \neq \varnothing$ since $x \in T_0$ and $y \in T_1$. If $z \in T_0$, then $g(z)_j \in X_{x_d}$  for every $j < d$, and if $z \in T_1$, then $g(z)_j \not\in X_{x_d}$  for every $j < d$.
Thus, the partition $\{g[T_0], g[T_1]\}$ of $g[T]$ demonstrates that $g[T]$ is not connected.
Thus, $g[T]$ is not an edge of  $L(X^d,P)$. 

In the second case, let $a \in X$ be such that $x_d = a$ for all $x \in T$. 
Thus, $g(x) \in X_a$ for each $x \in T$. 
Then, $g[T]$ is not an edge of $L(X_a^d,P)$, so it is not an edge of $L(X^d,P)$.

\smallskip

 {\em $P$ is not connected.} Let $P_0,P_1$ partition $P$ into two sets that 
demonstrate that $P$ is not connected. Thus, whenever $x \in P_0,y\in P_1$ and $j < d$, then 
$x_j \neq y_j$. Let $k_0 = |P_0|$ and $k_1 = |P_1|$. Thus, $k_0,k_1 < k$.

 By the inductive hypothesis, there are a $(d+1)$-dimensional $k_0$-template $Q_0$ and 
 a $(d+1)$-dimensional $k_1$-template $Q_1$ such that $e(Q_0) = e(P_0)$, $e(Q_1)= e(P_1)
$, $P_0 = \{x \harp d : x \in Q_0\}$,  $P_1 = \{x \harp d : x \in Q_1\}$ such that 
 $g$ is an embedding of $L(X^{d+1},Q_0)$ into $L(X^d,P_0)$ and also
 of $L(X^{d+1},Q_1)$ into $L(X^d,P_1)$. Moreover, we can arrange so that if $x \in Q_0$ and $y \in Q_1$, then $x_d \neq y_d$. 
 
 Let $Q = Q_0 \cup Q_1$. Clearly, $Q$ is a $(d+1)$-dimensional $k$-template 
 such that $e(Q) = e(P)$ and $P = \{x \harp d : x \in Q\}$. The ``moreover'' part of $(\bigstar_2)$ is easily checked. It remains to show that $g$ embeds $L(X^{d+1},Q)$ into $L(X^d,P)$. 
 Consider $T \in [X^d]^k$. Either $T$ is an edge of $L(X^{d+1},Q)$ or not.
 
 Suppose that $T$ is an edge of $L(X^{d+1},Q)$. Then, there is a homomorphism 
 $\varphi$ from $Q$ onto $T$. Then, $T_0 = \varphi[Q_0]$ and $T_1 = \varphi[Q_1]$ are homomorphic images of $Q_0$ and $Q_1$, respectively, so that $g[Q_0]$ and $g[Q_1]$ 
 are edges of $L(X^d,P_0)$ and $L(X^d,P_1)$. Making use of homomorphisms from $P_0$ and 
 $P_1$ onto $g[Q_0]$ and $g[Q_1]$, respectively, we can get a homomorphism from 
 $P$ onto $P$ onto $g[Q]$. 
 
 Suppose that $T$ is not an edge of $L(X^{d+1},Q)$. For a contradiction, assume that 
 $g[T]$ is an edge of $L(X^d,P)$, so let $\varphi : P \into g[T]$ be a surjective homomorphism. 
 Let $T_0 = g^{-1}\varphi(P_0)$ and $T_1 = g^{-1}\varphi(P_1)$. Without loss of generality, we can assume that $T_0$ is not an edge of $L(X^{d+1},Q_0)$. But then 
 $g[T_0]$ is not an edge of $L(X^d,P_0)$, implying that $g[T]$ is not an edge of $L(X^d,P)$, a contradiction. \qed
  
 \bigskip


 {\bf \S2.\@ Algebraic Hypergraphs.}   The infinite chromatic numbers of algebraic hypergraphs are determined in Theorem~2.2 of this section. This result leads to  a forbidden subhypergraph characterization in Corollary~2.4 of those algebraic $k$-hypergraphs $H$ for which $\chi(H) \leq \kappa$.

    The next definition is from \cite{avoid}.
 If  $A = A_0 \times A_1 \times \cdots \times A_{m-1}$, then a function $g : A \into Y$ is 
{\em one-to-one in each coordinate} if 
whenever $i < m$ and $a,b \in A$ are such that $a_j = b_j$ iff $i \neq j < m$, then 
$g(a) \neq g(b)$.

A subset $B \subseteq \RR^m$ is an {\it open box} if there are nonempty open intervals 
$B_0,B_1, \ldots, B_{m-1} \subseteq \RR$ such that $B = B_0 \times B_1 \times \cdots \times B_{m-1}$. Suppose that $P$ is an $m$-dimensional $k$-template, $B \subseteq \RR^m$ is an open box and 
 $H = (\RR^n,E)$ is a $k$-hypergraph. 
A function $f : B \into \RR^n$ is an {\it immersion} of $L(B,P)$ into $H$ 
if $f $ is a  semialgebraic analytic function that is  one-to-one in each coordinate and is such that 
whenever $\{x_0,x_1, \ldots, x_{k-1}\}$ is an edge of $L(B,P)$ and 
$f(x_0),f(x_1), \ldots, f(x_{k-1})$ are pairwise distinct, then $\{f(x_0),f(x_1), \ldots, f(x_{k-1})\}$ 
is an edge of $H$. If there is an immersion of $L(B,P)$ into $H$, then $L(B,P)$ is {\it immersible} in $H$. Observe that if $B \subseteq \RR^m$ is an open box, then 
$L(B,P)$ is immersible in $H$ iff $L(\RR^m,P)$ is immersible in $H$.

 \bigskip

{\sc Lemma 2.1}: {\em Suppose that  $H = (\RR^n,E)$ is an algebraic  $k$-hypergraph,    $P$ is a  $d$-dimensional $k$-template and  $L(\RR^d,P)$ is immersible in $H$. Then, $H$ contains an $L(\RR^d,P)$.}

\bigskip

{\it Proof}. Let $p(x_0,x_1, \ldots, x_{k-1})$ be a $(k,n)$-ary polynomial whose zero hypergraph is $H$, 
and let $f$ be an   immersion of  $L(\RR^d,P)$  into $H$. (This proof is easily modified so as not to make use of the analyticity of $f$.)

Let $T$ be a transcendence basis for $\RR$, and let $S \subseteq T$ be a finite set such that 
 $f$ is $S$-definable. Let $\{T_0,T_1, \ldots, T_{d-1}\}$ be a partition of $T \backslash S$ into $d$ sets each of which has cardinality $2^{\aleph_0}$. Then, $L(T_0 \times T_1 \times \cdots \times T_{d-1},P) \cong 
L(\RR^d, P)$. We claim that $f \harp(T_0 \times T_1 \times \cdots \times T_{d-1})$ embeds 
$L(T_0\times T_1 \times \cdots \times T_{d-1},P)$ into $H$. It suffices to show that $f$ is one-to-one 
on $T_0 \times T_1 \times \cdots \times T_{d-1}$. Suppose not, and let $\langle t_0,t_1, \ldots, t_{d-1} \rangle, \langle t_0',t_1', \ldots, t'_{d-1}\rangle \in T_0 \times T_1 \times \cdots \times T_{d-1}$ be distinct such that $f(t_0,t_1, \ldots, t_{d-1}) = f( t_0',t_1', \ldots, t'_{d-1})$. 
Without loss, suppose that $t_0 \neq t_0'$. Let $g : \RR \into \RR^n$ be such that 
$g(x) = f(x,t_1,t_2, \ldots,t_{d-1})$. Since $t_0$ is not algebraic over 
$S \cup \{t_1,t_2, \ldots, t_{d-1}\} \cup \{t'_0,t'_1, \ldots, t'_{d-1}\}$, there is an infinite set $X \subseteq \RR$ such that $g(x) = f(t'_0,t'_1, t'_2,\ldots, t'_{d-1})$ for every $x \in X$. But since $g$ is analytic, then 
$g(x)$ is constant on $\RR$, contradicting that $f$ is one-to-one in the first coordinate.  \qed

\bigskip

{\sc Theorem 2.2}: {\em Suppose that  $H$ is an algebraic $k$-hypergraph   and $\kappa$ 
is an infinite cardinal. 
 The following are equivalent$:$

\begin{itemize}

\item[$(1)$] $\chi(H) \leq \kappa;$

\item[$(2)$] whenever $P$ is a  $d$-dimensional $k$-template and  $H$ contains an $L(\RR^d,P)$, then  $\chi(L(\RR^d,P)) \leq \kappa;$

\item[$(3)$] whenever  $P$ is a  $d$-dimensional $k$-template and $L(\RR^d,P)$ is  immersible in $H$, then  $\chi(L(\RR^d,P)) \leq \kappa.$

\end{itemize}}

\bigskip

{\it Proof}. The implication $(1) \Longrightarrow (2)$ is trivial. The implication $(2) \Longrightarrow (3)$ easily follows from Lemma~2.1. 
We will  prove $(3) \Longrightarrow (1)$. The proof uses a technique from \cite{avoid}.

Since $(1)$  trivially holds when $\kappa \geq 2^{\aleph_0}$, we can assume that $\kappa < 2^{\aleph_0}$, 
although what follows does not depend on this inequality. 
Let $H = (\RR^n,E)$ and let $p(x_0,x_1, \ldots, x_{k-1})$ be a $(k,n)$-ary polynomial of which 
 $H$ is the zero hypergraph. Without loss of generality, we assume that $p(x_0,x_1, \ldots, x_{k-1})$ is 
 {\it symmetric}: whenever $a_0,a_1, \ldots, a_{k-1} \in \RR^n$ and 
 $\pi : k \into k$ is a permutation, then 
 $$p(a_0,a_1, \ldots, a_{k-1}) = 
  p(a_{\pi(0)},a_{\pi(1)}, \ldots, a_{\pi(k-1)}).$$
  This is possible since if  $p(x_0,x_1, \ldots, x_{k-1})$ is not symmetric, replace it with the symmetric  
  $\prod_{\pi} p(x_{\pi(0)},x_{\pi(1)}, \ldots, x_{\pi(k-1)})$. We also assume that $p(x_0,x_1, \ldots, 
  x_{k-1})$ is {\it reflexive}: whenever $a_0,a_1, \ldots, a_{k-1} \in \RR^n$ and 
  $|\{a_0,a_1, \ldots, a_{k-1}\}| < k$, then $p(a_0,a_1, \ldots, a_{k-1}) = 0$. 
   This is possible since if  $p(x_0,x_1, \ldots, x_{k-1})$ is not reflexive, replace it with the 
   reflexive  $p(x_0,x_1, \ldots, x_{k-1}) \prod_{i<j<k}(x_i - x_j)$. The assumptions of symmetry and reflexivity may not be essential, but they certainly don't hurt.

For each $ d < \omega$, let 
 $\theta_d : \RR^d \into \kappa$ be such that whenever $P$ is a $d$-dimensional $k$-template and $\chi(L(\RR^d,P)) \leq \kappa$, then $\theta_d$ is a proper coloring of $L(\RR^d,P)$. Such a $\theta_d$ exists since there are only finitely many non isomorphic $d$-dimensional $k$-templates. 
 Let $\FF \subseteq \RR$ be a countable real-closed subfield of $\RR$  such that all  coefficients 
  of $p$ are in $\FF$.

  Let $T$ be a transcendence basis for $\RR$ over $\FF$. 
 Then $|T| = 2^{\aleph_0}$. 
 For each $a = \langle a_0,a_1, \ldots, a_{n-1} \rangle \in \RR^n$, let $\supp(a)$,  the {\it support} of $a$, be the unique smallest 
 $S \subseteq T$ such that $\{a_0,a_1, \ldots, a_{n-1}\}$ is $(\FF \cup S)$-definable 
 (or, equivalently, such that each $a_i$ is algebraic over $\FF \cup S$).

We define a function $\varphi$ on $\RR^n$ as follows. Consider $a \in \RR^n$. 
Let $\supp(a) = \{t_0,t_1, \ldots, t_{d-1}\}$, where $t_0 < t_1 < \cdots < t_{d-1}$.
Then there are $q_0,q_1, \ldots, q_{d-1} , r_0,r_1, \ldots, r_{d-1} \in \QQ$  such that 
$$
q_0 < t_0 < r_0 < q_1 < t_1 < r_1 <  \cdots < q_{d-1} < t_{d-1} < r_{d-1}
$$
and an $\FF$-definable analytic  function 
$$
f:(q_0,r_0) \times (q_1, r_1) \times \cdots \times (q_{d-1},r_{d-1}) \into \RR^n
$$
such that  $f$ is one-to-one in each coordinate and $f(t_0,t_1, \ldots,t_{d-1}) = a$. 
(This $f$ is the {\em determining} function for $a$). Observe that there are only countably many 
determining functions since each one is $\FF$-definable and $\FF$ is countable. 
Finally, we let 
$\varphi(a) = \big\langle f, \theta_d( t_0,t_1, \ldots, t_{d-1} )  \big\rangle$. 

\smallskip

Clearly, $\varphi$ is a $\kappa$-coloring of $H$. 
We now claim that it is a proper coloring. For a contradiction, suppose that 
$\{a_0,a_1, \ldots, a_{k-1}\} \in E$ and 
$\varphi(a_0) = \varphi(a_1) = \cdots = \varphi(a_{k-1}) = \langle f, \alpha \rangle$.

For each $i < k$, let $t_i = \langle t_{i,0}, t_{i,1}, \ldots, t_{i,d-1} \rangle \in \dom(f)$ be such 
that $\supp(a_i) = \{ t_{i,0}, t_{i,1}, \ldots, t_{i,d-1}\}$. Thus, $a_i = f(t_{i,0}, t_{i,1}, \ldots, t_{i,d-1})$.
Let  $P = \{ t_0,t_1, \ldots, t_{k-1} \}$. Since the $a_i$'s are  pairwise distinct, so are the $t_i$'s. Hence, $P$ is a $d$-dimensional $k$-template. 

Since $t_0,t_1, \ldots, t_{k-1} \in \RR^d$ are distinct, $\theta_d(t_0) = \theta_d(t_1) = \cdots = \theta_d(t_{k-1}) = \alpha$ 
and $\{ t_{0}, t_{1}, \ldots, t_{k-1} \}$ is an edge of $L(\RR^d,P)$, then 
$\theta_d$ is not a proper coloring of $L(\RR^n,P)$. 
 Thus, 
 it must be that $\chi(L(\RR^d,P)) > \kappa$.

We will arrive at a contradiction by proving that $f$ is an immersion of $L(\dom(f),P)$ is immersible into $H$. 
 This is a consequence of the following claim.

\smallskip

{\sf Claim}: {\em If $\{s_0,s_1, \ldots, s_{k-1}\} \subseteq \dom(f)$ is a 
$d$-dimensional $k$-template that is a homomorphic image of $P$ and $f(s_0),f(s_1), \ldots, f(s_{k-1})$ are pairwise distinct, then   $\{ f(s_0),f(s_1), \ldots, f(s_{k-1}) \} \in E$}. 

\smallskip

To prove the claim, suppose that $s_0,s_1, \ldots, s_{k-1} \in \dom(f)$, $\{s_0,s_1, $ $\ldots, s_{k-1}\}$ is a 
$k$-template and the function 
 $t_i \mapsto s_i$   demonstrates that $\{s_0,s_1, $ $\ldots, s_{k-1}\}$ is a 
homomorphic image of $P$. 
We have that  $$p(f(t_0),f(t_1), \ldots, f(t_{k-1})) = 0$$ from which it follows,  by \cite[Lemma~2.4]{avoid}, that $$p(f(s_0),f(s_1), \ldots, f(s_{k-1})) = 0.$$ Therefore,  $\{ f(s_0),f(s_1), \ldots, f(s_{k-1}) \} \in E$. \qed

\bigskip

{\sc Corollary 2.3}: {\em Suppose that $2 \leq k < \omega$ and $\kappa$ is a cardinal. 
 The following are equivalent$:$

$(1)$ there is an algebraic $k$-hypergraph $H$ such that $\chi(H) = \kappa;$

$(2)$ $ 1 \leq \kappa \leq \aleph_0$ or $\kappa \leq 2^{\aleph_0} \leq \kappa^{+(k-2)}$.} 

\bigskip

{\it Proof}. $(1) \Longrightarrow (2)$: Let $H$ be an algebraic $k$-hypergraph and let $\kappa =  \chi(H)$. Obviously, $\kappa \leq 2^{\aleph_0}$. We can assume that $\aleph_0< \kappa < 2^{\aleph_0}$ as otherwise (2) is true. Since $\chi(H) > \aleph_0$, it follows from $(2) \Longrightarrow (1)$ of Theorem~2.2 that there   are $d < \omega$ and a  
$d$-dimensional $k$-template $P$ such that $H$ contains an $L(\RR^d,P)$ and $\chi(L(\RR^d,P)) > \aleph_0$. By Theorem~1.1, $\chi(L(\RR^d,P))^{+(e(P)-1)} \geq 2^{\aleph_0}$. Since $\kappa \geq 
\chi(L(\RR^d,P))$ and $e(P) \leq k-1$, then $\kappa^{+(k-2)} \geq 2^{\aleph_0}$.

\smallskip

$(2) \Longrightarrow (1)$: If $1 \leq \kappa < \aleph_0$, then let $H' = (V,E)$ be a finite 
$k$-hypergraph such that $V \subseteq \RR$ and $\chi(H') = \kappa$. (For instance, let $H'$ 
be a complete $k$-hypergraph with $|V| = (\kappa-1)(k-1) +1$.) Then
$H = (\RR,E)$ is an algebraic $k$-hypergraph and $\chi(H) = \kappa$.

If $\kappa \leq 2^{\aleph_0} \leq \kappa^{+(k-2)}$, then  Corollary~1.6 implies that there is an algebraic $k$-hypergraph $H$ such that  $\chi(H) = \kappa$. 

Next, suppose that $\kappa = \aleph_0$ and $k=2$. Let $p(x_0,x_1,y_0,y_1)$ be the $4$-ary polynomial $x_1-y_0$. 
Consider it as a $(2,2)$-ary polynomial $p(x,y)$, where $x = \langle x_0,x_1 \rangle$ 
and $y = \langle y_0,y_1 \rangle$, and let $H = (\RR^2,E)$ be its zero graph. 
Then, $H$ is an algebraic graph.  We first show that $\chi(H) \leq \aleph_0$ by exhibiting a  proper 
$\aleph_0$-coloring $\varphi$ of it. Let $Q_0, Q_1$ be two disjoint dense sets of nonzero rationals. Then let  $\varphi : \RR^2 \into \QQ$ be such that whenever  $a<b\in \RR$, 
then $\varphi(a,b) \in Q_0$, $\varphi(b,a) \in Q_1$, $a < \varphi(a,b), \varphi(b,a) < b$ and $\varphi(a,a) = 0$.  It is clear that $\varphi$ is a proper $\aleph_0$-coloring.
To see that $\chi(H) \geq \aleph_0$, consider any coloring $\psi : \RR^2 \into m < \omega$. 
By Ramsey's Theorem, there are integers $a < b < c$ such that $\psi(a,b) = \psi(b,c)$, so $\psi$ is not proper.

For $k \geq 3$, it is trivial to obtain from the previous graph $H$ an algebraic $k$-hypergraph whose chromatic 
number is~$\aleph_0$. \qed

\bigskip

The next corollary shows that  those algebraic 
$k$-hypergraphs that are $\kappa$-colorable, where  $\kappa \geq \aleph_0$, can be characterized by a finite set of forbidden 
algebraic subhypergraphs.

\bigskip

{\sc Corollary 2.4}: {\em Suppose that $2 \leq k < \omega$ and $ \kappa$ is an infinite cardinal. There is a finite set ${\mathcal F}$ of 
algebraic $k$-hypergraphs such that  for any algebraic $k$-hypergraph $H$, 
$\chi(H) \leq \kappa$ iff $H$ does not contain any $F \in {\mathcal F}$.}

\bigskip

{\it Proof}. Let ${\mathcal F}$ be the set of all $L(\RR^e,P)$, where $e \leq k-1$ and $P$ is an
$e$-dimensional $k$-template such that $\chi(L(\RR^e,P)) > \kappa$. Clearly, ${\mathcal F}$ is a finite set of algebraic $k$-hypergraphs. 

Let $H$ be an algebraic $k$-hypergraph. If   $\chi(H) \leq \kappa$ and  $F \in {\mathcal F}$, then $H$ does not contain an $F$ since $\chi(F) > \kappa \geq \chi(H)$. Conversely, suppose that $\chi(H) > \kappa$. 
By $(2) \Longrightarrow (1)$ of Theorem~2.2, there is a $d$-dimensional $k$-template $Q$ 
such that $H$ contains $L(\RR^d,Q)$ and $\chi(L(\RR^d,Q)) > \kappa$. Let $e = e(Q)$. 
Clearly, $e \leq k-1$. 
By Theorem~1.1, $\kappa^{+(e-1)} < 2^{\aleph_0}$. Following Lemma~1.9, let $P$ be an $e$-dimensional $k$-template  such that $L(\RR^d,Q)$ contains an $L(\RR^e,P)$.  Let $F = L(\RR^e,P)$. Then $F$ is embeddable into $H$ since it is embeddable into $L(\RR^d,Q)$. Since \mbox{$e(P) \leq e$},  Theorem~1.1 implies that 
$\chi(F) > \kappa$. Thus, $F \in {\mathcal F}$. \qed

\bigskip

{\sc Remark 2.5:} The proof of Corollary~2.4 shows that all  hypergraphs in ${\mathcal F}$ can have the form 
$L(\RR^e,P)$, where $P$ is an $e$-dimensional $k$-template and   $e \leq k-1$. Moreover,  we can also require that  $P$ be simple. By Lemma~1.10, we can also arrange for ${\mathcal F}$ to consist only of  hypergraphs having the form  $L(\RR^{k-1},P)$, where $P$ is a $(k-1)$-dimensional $k$-template.

\bigskip

The next corollary asserts that master colorings exist. 
 Komj\'ath \cite{kom92} first considered master colorings and proved,  assuming  {\sf CH}, that they exist for   $\kappa = \aleph_0$. The {\sf CH} assumption was subsequently eliminated in~\cite{avoid}.
 
\bigskip

{\sc Corollary 2.6}: {\em Suppose  that $1 \leq n < \omega$ and $\kappa$ is an infinite cardinal.
 For each $\alpha < \kappa$, let  
$H_\alpha$ be an algebraic hypergraph on $\RR^n$ such that $\chi(H_\alpha) \leq \kappa$. 
Then there is a function $\varphi: \RR^n \into \kappa$ that is a proper coloring of $H_\alpha$ 
 for each $\alpha < \kappa$,} 

\bigskip

{\it Proof}. We will use not only Theorem~2.2 but also its proof. 

We can assume that $\kappa < 2^{\aleph_0}$. Suppose, for each $\alpha < \kappa$, that  
 $H_\alpha = (\RR^n,E_\alpha)$ is the zero hypergraph of  the $(k_\alpha,n)$-polynomial $p_\alpha(x_0,x_1, \ldots,$ $x_{k_\alpha-1})$. 

Observe that the field $\FF$ 
in the proof of Theorem~2.2 was required to be countable, but it would have done no harm if it had been of  cardinality $\kappa$. 
So, in that proof, choose $\FF$ so that $|\FF| = \kappa$ and  all coefficients of 
 each $p_\alpha$ are in $\FF$. Next, observe that the definition of the $\kappa$-coloring 
 $\varphi$ did not depend on $H$. Thus, the same $\varphi$ is a proper $\kappa$-coloring 
 for every $k$-hypergraph $H$ for which $\chi(H) \leq \kappa$ and $H$ is the zero hypergraph of  
  a $(k,n)$-ary  polynomial over $\FF$. In particular, $\varphi$ is a proper $\kappa$-coloring of each~$H_\alpha$. \qed
 
 \bigskip
 
 There is an instance of Corollary~2.6 for graphs  that merits special mention. 
  If $n < \omega$ and $D$ is a set of positive reals, then the $D$-distance graph on $\RR^n$,
 which generalizes the unit distance graph,  is the graph 
${\mathbf X}_n(D)$ whose vertices are the points in $\RR^n$ and whose edges are those pairs of points at a distance in $D$.
This graph is algebraic iff $D$  is either finite or  the set of all positive reals. Komj\'ath \cite{kom} showed  that $\chi({\mathbf X}_n(D)) \leq \aleph_0$ whenever $D$ is countable. The next corollary   extends  Komj\'{a}th's result to  arbitrary ~$D$.

 \bigskip

{\sc Corollary~2.7}: {\em If $D$ is a  set of positive reals, then $\chi({\mathbf X}_n(D)) \leq |D| + \aleph_0$.} \qed

\bigskip

 It easily follows that if $D$ is the positive part of  some nontrivial additive subgroup of $\RR$,  then $\chi({\mathbf X}_n(D)) = |D|$. In contrast to this, it was proved in 
\cite{ks} that if $D$ is an algebraically independent set of positive reals, then  $\chi({\mathbf X}_n(D)) \leq \aleph_0$.  Bukh \cite{bukh} made the conjecture (still open for $n \geq 2$) that 
$\chi({\mathbf X}_n(D)) < \aleph_0$ for such $D$. Getting some $D$ such that $|D| = 2^{\aleph_0}$ 
and $\chi({\mathbf X}_n(D)) < \aleph_0$ is an easy matter: just let $D = [a,b]$, where $0 < a < b$. 

\bigskip

There is  a metatheorem that can be very roughly stated as: If an 
 $(k,n)$-ary polynomial $p(x_0,x_1, \ldots, x_{k-1})$ and $m < \omega$ are such that the sentence

\smallskip
 
\noindent $(*)$ \hspace{40pt} {\em $p(x_0,x_1, \ldots, x_{k-1})$ is avoidable iff $2^{\aleph_0} \leq \aleph_m$}

\smallskip

\noindent is provable, then the seemingly stronger sentence

\smallskip
 
\noindent $(**)$ \hspace{18pt} {\em For any $\kappa$, $p(x_0,x_1, \ldots, x_{k-1})$ is $\kappa$-avoidable iff $\kappa^{+m} \geq 2^{\aleph_0}$}

\smallskip
\noindent is also provable. As an illustrative example, we show how the following corollary can be obtained from Fox's result mentioned in the introduction. (To be fair, it should be noted that this corollary is not new, being implicit in \cite{fox}.)

\bigskip

{\sc Corollary 2.8:} {\em If $k < \omega$, then the $(k+3,1)$-ary polynomial $x_0 + x_1 + \cdots + x_k - x_{k+1} - kx_{k+2}$ is $\kappa$-avoidable iff $\kappa^{+k} \geq 2^{\aleph_0}$.}

\bigskip

{\it Proof}. (In this proof we will refer to models of ${\sf ZFC}$, but what we really mean are models of some explicitly given finite fragment of ${\sf ZFC}$ that is large enough to prove all 
the relevant facts that are needed.)  Fix $k < \omega$. Let $p(x)$ be the given $(k+3)$-ary polynomial and let $H$ be its zero hypergraph.  It suffices to show that the sentence

\smallskip

\noindent (1) \hspace{30pt}  {\em for any cardinal $\kappa$, $H$ is $\kappa$-colorable iff $\kappa^{+k} \geq 2^{\aleph_0}$}

\smallskip

\noindent is true in every countable model of ${\sf ZFC}$. Suppose that $M$ is such a model. 
Let $M[G]$ be a generic extension in which $2^{\aleph_0} = \aleph_{k+1}$. According to Fox, 
$p(x)$ is not avoidable in $M[G]$. Then, from Theorems~2.2 (with $\kappa = \aleph_0$) and 1.1, we get that the sentence 

\smallskip

\noindent (2) \hspace{15pt}{\em there is a $d$-dimensional $(k+3)$-template $P$ such that}  

{\hspace {18pt} {\em $L(\RR^d,P)$ is immersible in 
$H$ and  $e(P) = k+2$, but}

{\hspace{18pt} {\em none with $e(P) = k+1$}

\smallskip

\noindent is true in $M[G]$. By  Tarski's Theorem on the decidability of Th$(\widetilde \RR)$, sentence (2) is equivalent to an arithmetic sentence (in fact, a $\Delta^0_2$ sentence), so it is absolute. Hence,  it is also true in $M$. Then, again using 
Theorems~2.2 and 1.1, we conclude that $(1)$ is true in $M$. \qed 

\bigskip

We end by stating without proof  one other consequence of the metatheorem based on \cite[Prop.\@ 1.4]{avoid}.

\bigskip

{\sc Corollary 2.9:} {\em For $2 \leq n < \omega$, let $H = (\RR^n,E)$ be the $(n+1)$-hypergraph 
whose edges consist exactly  of those $A \subseteq [\RR^n]^{n+1}$ that are the vertices of some orthogonal $n$-simplex. Then, for any cardinal $\kappa$, $H$ is $\kappa$-colorable iff 
$\kappa^{+(n-1)} \geq 2^{\aleph_0}$.} \qed

\bibliographystyle{plain}

\end{document}